\documentclass[final,5p,twocolumn,times,numbers,sort&compress]{elsarticle}

\usepackage{amssymb}
\usepackage[fleqn]{amsmath}
\usepackage{algorithm}
\usepackage{algpseudocode}
\usepackage{graphicx}
\usepackage{booktabs}
\usepackage{siunitx}
\usepackage{hyperref}
\usepackage{float}
\usepackage{enumitem}
\usepackage{multirow}
\usepackage{multicol} 
\usepackage{nomencl}
\usepackage{framed}
\makenomenclature 
\renewcommand*\nompreamble{\begin{multicols}{2}}
\renewcommand*\nompostamble{\end{multicols}}
\renewcommand{\nomgroup}[1]{
\ifthenelse{\equal{G}{#1}}{\item[\textit{Greek symbols}]}{}
\ifthenelse{\equal{S}{#1}}{\item[\textit{Subscripts}]}{}
\ifthenelse{\equal{U}{#1}}{\item[\textit{Superscripts}]}{}}

\setlist{leftmargin=15pt,labelindent=15pt}
\setlist[enumerate]{wide=0pt, leftmargin=15pt, labelwidth=15pt, align=left}

\hypersetup{
  colorlinks   = true, 
  urlcolor     = blue, 
  linkcolor    = blue, 
  citecolor   = blue 
}

\journal{Elsevier}

\begin{document}

\begin{frontmatter}



\title{Topology optimization for microchannel heat sinks with nanofluids using an Eulerian-Eulerian approach}


\author[a]{Chih-Hsiang Chen}
\author[a]{Kentaro Yaji\corref{c1}}
\cortext[c1]{Corresponding author}
\ead{yaji@mech.eng.osaka-u.ac.jp}
\affiliation[a]{organization={Department of Mechanical Engineering, Graduate School of Engineering, Osaka University},
            addressline={2-1, Yamadaoka}, 
            city={Suita},
	   state={Osaka},
            postcode={ 565-0871}, 
            country={Japan}}

\begin{abstract}
The demand for high-performance heat sinks has significantly increased with advancements in computing power and the miniaturization of electronic devices. Among the promising solutions, nanofluids have attracted considerable attention due to their superior thermal conductivity. However, designing a flow field that effectively utilizes nanofluids remains a significant challenge due to the complex interactions between fluid and nanoparticles. In this study, we propose a density-based topology optimization method for microchannel heat sink design using nanofluids. An Eulerian–Eulerian framework is utilized to simulate the behavior of nanofluids, and the optimization problem aims to maximize heat transfer performance under a fixed pressure drop. In numerical examples, we investigate the dependence of the optimized configuration on various parameters and apply the method to the design of a manifold microchannel heat sink. The parametric study reveals that the number of flow branches increases with the increased pressure drop but decreases as the particle volume fraction increases. In the heat sink design, the topology-optimized flow field achieves an 11.6$\%$ improvement in heat transfer performance compared to a conventional parallel flow field under identical nanofluid conditions.
\end{abstract}


\end{frontmatter}

\nomenclature{$\bold{u}$}{velocity vector ($\rm m/s$)}
\nomenclature{$u$}{velocity component in the $x$-direction ($\rm m/s$)}
\nomenclature{$v$}{velocity component in the $y$-direction ($\rm m/s$)}
\nomenclature{$p$}{pressure (Pa)}
\nomenclature{$T$}{temperature (K)}
\nomenclature{$c$}{specific heat capacity ($\rm J / kg \cdot K$)}
\nomenclature{$k$}{thermal conductivity ($\rm W / m \cdot K$)}
\nomenclature{$K$}{momentum exchange coefficient ($ \rm kg / m^{3} \cdot s$ )}
\nomenclature{$h$}{volumetric heat exchange coefficient ($\rm W / m^{3} \cdot K$ )}
\nomenclature{$h_{\rm p}$}{fluid-particle heat transfer coefficient ($\rm W / m^{2} \cdot K$ )}
\nomenclature{$Q$}{volumetric heat generation coefficient ($\rm W / m^{3} \cdot K$ )}
\nomenclature{$C_{\rm d}$}{drag coefficient}
\nomenclature{$d_{\rm p}$}{particle diameter (m)}
\nomenclature{$Nu$}{Nusselt number}
\nomenclature{$Re$}{fluid Reynolds number}
\nomenclature{$Re_{\rm p}$}{particle Reynolds number}
\nomenclature{$Pr$}{Prandtl number}
\nomenclature{$L$}{inlet length (m)}
\nomenclature{$H$}{plate length (m)}
\nomenclature{$\bold{x}$}{position vector (m)}
\nomenclature{$R$}{filter radius (m)}
\nomenclature{$q$}{convex parameter}
\nomenclature{$J$}{objective function}
\nomenclature{$D$}{design domain}
\nomenclature{$i$}{number of flow simulation iterations}
\nomenclature{$n$}{number of optimization iterations}
\nomenclature{tol}{convergence tolerance for flow simulation}
\nomenclature[G]{$\mu$}{viscosity ($\rm Pa \cdot s$)}
\nomenclature[G]{$\rho$}{density ($\rm kg / m^{3}$)}
\nomenclature[G]{$\phi$}{volume fraction}
\nomenclature[G]{$\alpha$}{inverse permeability ($\rm kg/m^{3} \cdot s$)}
\nomenclature[G]{$\gamma$}{design variable}
\nomenclature[G]{$\Omega$}{fluid domain}
\nomenclature[G]{$\epsilon$}{convergence tolerance for optimization}
\nomenclature[G]{$\beta$}{projection parameter}
\nomenclature[G]{$\eta$}{threshold parameter}
\nomenclature[G]{$\Gamma$}{boundary of the analysis domain}
\nomenclature[S]{f}{fluid phase}
\nomenclature[S]{p}{particle phase}
\nomenclature[S]{ave}{average}
\nomenclature[U]{in}{inlet}
\nomenclature[U]{out}{outlet}
\nomenclature[U]{w}{wall}
\begin{table*}[!t]
  \begin{framed}
    \printnomenclature
  \end{framed}
\end{table*}

\section{Introduction}

The increasing computational demands and trend toward highly integrated circuits have led to a rise in the power consumption of electronic devices. Thermal dissipation has become a critical concern due to high heat flux generation and the miniaturization of electronic devices~\cite{he2021thermal}. To address this issue, the microchannel heat sink (MCHS), first proposed by Tuckerman and Pease~\cite{tuckerman1981high}, has become an attractive solution due to its high heat transfer performance and compact size~\cite{gao2022fluid}. Over the past decades, numerous studies have sought to improve MCHS from various perspectives, including operating conditions~\cite{rahman2000measurements, qu2002experimental}, coolant types~\cite{maghrabie2023microchannel, saeed2021numerical}, material selection~\cite{kepekci2020comparative, muhammad2020comparison, kumar2018material}, and channel configurations~\cite{kose2022parametric, pan2020study, dash2022role}.

In terms of coolants, nanofluids introduced by Choi and Eastman~\cite{choi1995enhancing} have attracted significant attention due to their superior heat transfer properties. The fundamental concept of nanofluids is to disperse nanoscale particles into a base fluid to enhance heat transfer performance, and they have been widely applied in the field of heat dissipation~\cite{chein2005analysis, jang2006cooling, sohel2014experimental}. Although typical flow field designs—such as parallel~\cite{ho2013experimental}, fin-and-pin~\cite{hasan2014investigation}, and wavy~\cite{sajid2019experimental}—are widely employed with nanofluids, the complex interactions between nanoparticles and the base fluid present significant challenges in designing efficient flow fields for nanofluids. Meanwhile, numerous studies have conducted parametric investigations based on these designs~\cite{lee2015fluid, wu2014parametric, chai2016parametric}; however, their limited design freedom has constrained the performance improvements in MCHS.

As a method with high design flexibility, topology optimization can generate optimized configurations without relying on the designer's intuition. In topology optimization, the density method employs a continuous design variable field, which enables smooth variations of material distribution throughout the design domain during the optimization process. The introduction of a penalization factor eliminates intermediate values of these design variables, resulting in a binary distribution of material~\cite{bendsoe2013topology}. Based on the idea of the density method, Borrvall and Petersson proposed a topology optimization approach for Stokes flow, introducing a fictitious force based on design variables as a no-slip boundary condition to determine the material distribution of the pseudo-solid region~\cite{borrvall2003topology}. Further research has also developed topology optimization methods for heat sink design. Dede~\cite{dede2009multiphysics} and Yoon~\cite{yoon2010topological} first proposed a topology optimization method for forced convection heat transfer problems. Koga et al. proposed a topology optimization method for heat sinks, defining a multi-objective function to maximize heat dissipation while simultaneously minimizing power dissipation in their study~\cite{koga2013development}. Matsumori et al. developed a topology optimization method for heat sinks that accounts for constrained input power~\cite{matsumori2013topology}. Yaji et al. developed a topology optimization method to maximize heat exchange, using a level set boundary expression to control the geometric complexity of the optimized configurations~\cite{yaji2016topology}. Dilgen et al. proposed a topology optimization method for turbulent forced convection heat transfer systems~\cite{dilgen2018density}. Zeng et al. developed a two-layer model for topology optimization to reduce the computational cost of 3D heat sink design and experimentally verified that the optimized result exhibits superior performance~\cite{zeng2018experimental}. Zhou et al. introduced a topology optimization method for manifold microchannel heat sinks and investigated the optimized results numerically and experimentally~\cite{zhou2024topologically}. Although these studies proposed topology optimization methods for heat sinks, they did not consider nanofluids as potential coolants. In terms of nanofluids, Zhang et al. proposed a topology optimization method for MCHS~\cite{zhang2020topology}. Due to the diluteness of nanoparticles, a single-phase model is employed to simulate the behavior of nanofluids, with its physical properties being temperature-dependent. In their study, while optimized configurations can be obtained using a single-phase model, the performance of the structures optimized with nanofluids was inferior to those optimized with pure water.

In the field of numerical investigations for nanofluids, although the single-phase model is frequently used, recent studies have increasingly focused on the two-phase model for nanofluids. In the two-phase model, the particle phase is treated as a continuous phase to account for the complex interactions between nanoparticles and the fluid. Bianco et al. investigated the laminar forced convection flow of $\rm Al_{2}O_{3}$ nanofluid in a circular pipe with the single-phase and two-phase models. The maximum difference in the average heat transfer coefficient between the single-phase model and the two-phase model is 11$\%$ in their study~\cite{bianco2009numerical}. Kalteh et al. conducted numerical simulations of nanofluid forced convection in a microchannel using the Eulerian-Eulerian approach. The average Nusselt number enhancement obtained with the Eulerian-Eulerian approach is higher than that achieved with the single-phase model~\cite{kalteh2011eulerian}. In terms of experimental validation, Fard et al. compared the simulation results of single-phase and two-phase models for nanofluids with copper nanoparticles and found that the two-phase model results aligned more closely with experimental data~\cite{fard2010numerical}. Beheshti et al. conducted a comparative numerical study of heat transfer with turbulent nanofluids in an annular channel. In their study, the single-phase model underestimates the experimental data at concentrations (1$\%$ and 1.5$\%$). In contrast, the two-phase model yields results that are closer to the experimental values~\cite{beheshti2015comparative}. In summary, these studies have demonstrated that the two-phase model provides more accurate results and serves as an effective alternative for simulating the heat transfer of nanofluids.

Several studies have employed topology optimization to address particle flow problems. In these studies, the Eulerian-Lagrangian approach has been extensively utilized to optimize particle trajectory and manipulation~\cite{andreasen2020framework, yoon2020transient, yoon2022transient}. However, this approach leads to high computational costs when applied to the collective behavior of nanofluids. In contrast, although Chen et al. demonstrated that the Eulerian–Eulerian approach can effectively tackle topology optimization problems for particle flows~\cite{chen2024topology}, the assumption of one-way coupling in their study limits the ability to capture the bidirectional interactions between the fluid and particles in nanofluids.

\begin{figure}[h]
\centering
\includegraphics[width=0.45\textwidth]{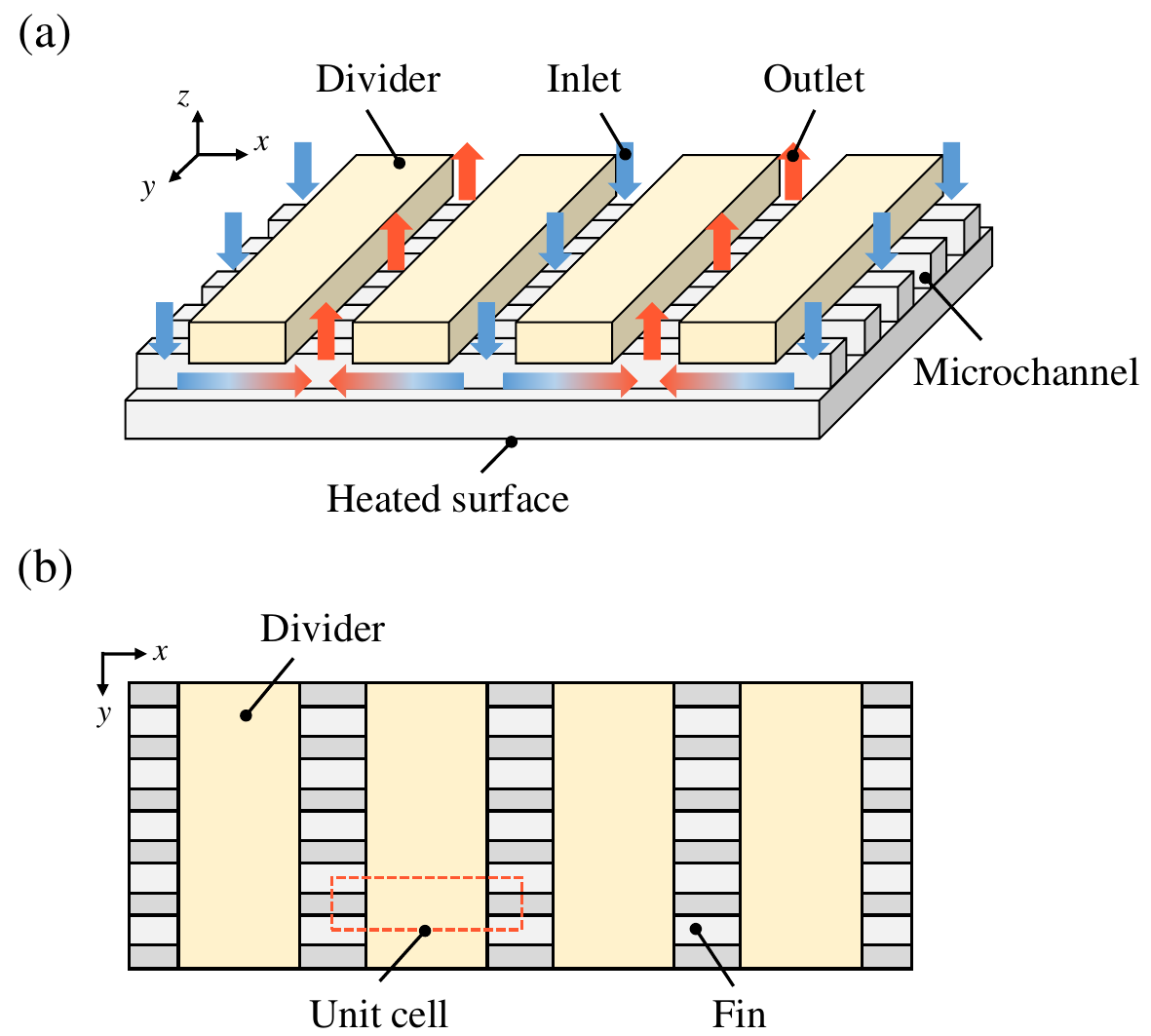} 
\caption{Schematic diagram of a manifold microchannel heat sink: (a) The flow direction and key components in the MMCHS (b) Top view of a unit cell in the MMCHS.}
\label{Fig.0}
\end{figure}

Based on the above, we propose a topology optimization method for flow field design in MCHS with nanofluids using an Eulerian-Eulerian approach. In the numerical model, the particle phase is treated as a continuous phase, and fluid-particle bidirectional interactions are incorporated into the momentum and energy equations for both phases based on the Eulerian-Eulerian approach. The optimization problem is defined as a heat transfer maximization problem under a fixed pressure drop. The sensitivity information is derived using automatic differentiation, and GCMMA (Globally Convergent Method of Moving Asymptotes) is employed for solving the optimization problem~\cite{svanberg2002class}. This study includes two numerical examples. In the first example, we optimize the flow field in a square design domain commonly used for topology optimization in thermal management. This example serves as a preliminary exploration of various parameters to validate the proposed optimization method. In the second numerical example, we consider a manifold microchannel heat sink (MMCHS) for optimization, as shown in Fig.~\ref{Fig.0}. In the MMCHS, dividers segment the microchannels into multiple sections, allowing fluid to flow downward into each section and upward out of it. Since the MMCHS consists of periodic unit cells, we optimize the flow field in a unit cell with a simplified two-dimensional model in this study.

The subsequent sections of this paper are organized as follows:  In Section~\ref{Analysis model}, we demonstrate the formulation of the governing equations and boundary conditions. In Section~\ref{Problem formulation}, we describe the representation of fluid and solid regions in the topology optimization method and formulate the optimization problem for flow field design in MCHS. In Section~\ref{Numerical implementation}, we provide an overview of the algorithm used in the fluid solver and outline the optimization process. In Section~\ref{Results and discussion}, we compare the effects of different parameters on the optimized flow field and optimize the flow field in the unit cell of MMCHS. Finally, we present the summary of findings and future work in Section~\ref{Conclusion}.
\section{Analysis model}
\label{Analysis model}
\subsection{Assumptions}
\label{Assumptions}
The particle volume fraction has a significant impact on the behavior of nanofluids. Although increasing the particle volume fraction can enhance heat transfer, excessive concentration may lead to nanoparticle aggregation, which increases viscosity and reduces heat transfer efficiency~\cite{chakraborty2020stability}. Therefore, the volume fraction of nanofluids is typically low to prevent these negative effects. In this study, we adopt the two-way coupling assumption, meaning that the behaviors of the fluid and particles influence each other, while interactions between particles are neglected. Additionally, we consider a steady-state laminar flow for the numerical examples due to the broad range of laminar flow applications in MCHS~\cite{lee2002experimental, chai2016numerical, zhu2022numerical}. Although the fluid in the MMCHS exhibits both upward and downward flow behaviors, under high aspect ratio conditions—where the channel width is significantly greater than its depth—the distributions of physical quantities remain sufficiently uniform~\cite{li2006aspect, kamal2017analysis}. Consequently, the flow field can be approximated as two-dimensional. In addition, we further consider a uniform heat distribution along the channel's height and simplify the volumetric heat generation into the two-dimensional energy equation. The basic assumptions are summarized as follows:

\begin{itemize}
\item The flow fields are treated as two-dimensional.
\item The fluid is treated as a steady incompressible laminar flow.
\item Heat distribution is uniform along the channel's height.
\item Particle-particle interactions are not considered.
\item The gravity effect is not considered.
\item Shape effects are neglected, and particles are spherical.
\item Particle aggregation is not considered.
\end{itemize}
\begin{figure}[h]
\centering
\includegraphics[width=0.35\textwidth]{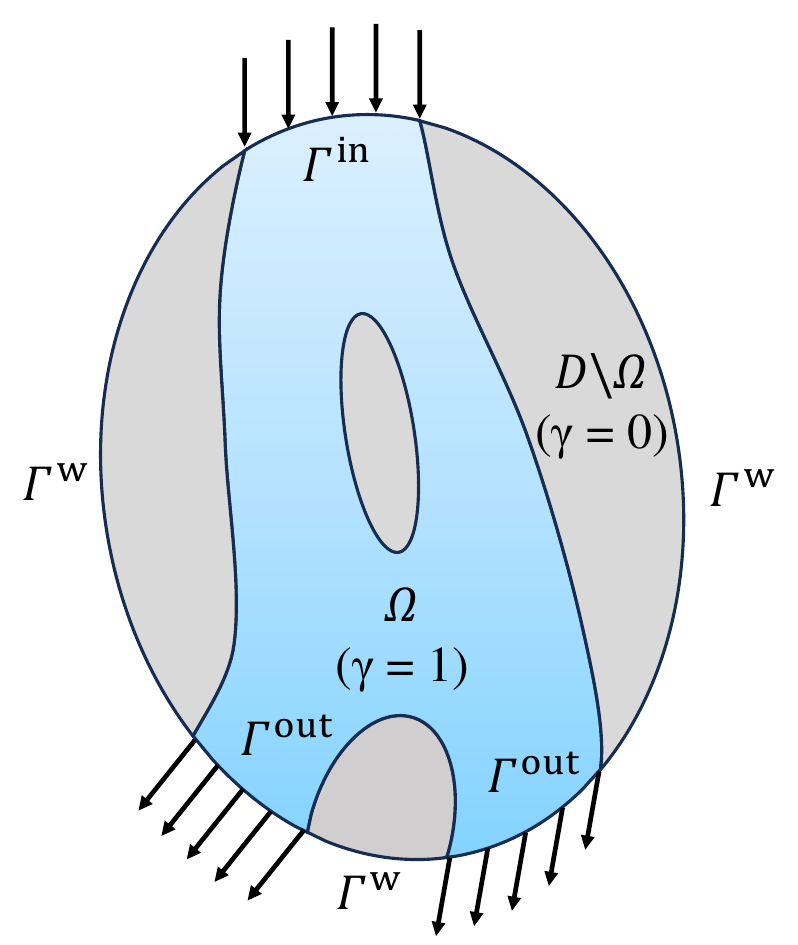} 
\caption{Schematic representation of the fluid domain $\Omega$ and the solid domain $D \backslash \Omega$ in the fixed design domain $D$, including the corresponding boundary expressions for the inlet $\mathit{\Gamma}^{\rm in}$, outlet $\mathit{\Gamma}^{\rm out}$, and wall $\mathit{\Gamma}^{\rm w}$.}
\label{Fig.1}
\end{figure}
\subsection{Governing equations}
\label{Governing equations}
Based on the assumptions in Section~\ref{Assumptions} and the Eulerian-Eulerian model proposed by Kalteh et al.~\cite{kalteh2011eulerian}, the governing equations for the fluid phase, including continuity, momentum, and energy, are as follows:
\begin{align}
\label{gef}
&\nabla \cdot (\phi_{\rm f} \bold{u}_{\rm f}) = 0,  \\[10pt]
\label{nsf}
&\nabla \cdot ({\phi_{\rm f} \rho_{\rm f} \bold{u}_{\rm f} \bold{u}_{\rm f}})  =  -\phi_{\rm f} \nabla p + \nabla \cdot (\mu_{\rm f} \phi_{\rm f}(\nabla \bold{u}_{\rm f} + \nabla \bold{u}^\intercal_{\rm f})) \\
& \ \ \ \ \ \ \ \ \ \ \ \ \ \ \ \ \ \ \ \ \ \ + K(\bold{u_{\rm p}} - \bold{u_{\rm f}}) \nonumber, \\
\label{geef}
&\nabla \cdot ({\phi_{\rm f} \rho_{\rm f}  c_{\rm f}\bold{u}_{\rm f} T_{\rm f}})  =\nabla \cdot (\phi_{\rm f}\hat{k}_{\rm f}\nabla T_{\rm f}) + h(T_{\rm p} - T_{\rm f}),
\end{align}
where $\bold{u}$ is the velocity vector, $p$ is the pressure, $\mu$ is the viscosity, $\rho$ is the density, $\phi$ is the volume fraction, $K$ is the momentum exchange coefficient, $T$ is the temperature, $c$ is the specific heat capacity, $h$ is the heat exchange coefficient, and $\hat{k}$ is the effective thermal conductivity. The subscripts `f' and `p' denote the fluid and particle phases. Similarly, the governing equations for the particle phase are described as follows:
\begin{align}
\label{gep}
&\nabla \cdot (\phi_{\rm p} \bold{u}_{\rm p}) = 0,  \\[10pt]
\label{nsp}
&\nabla \cdot ({\phi_{\rm p} \rho_{\rm p} \bold{u}_{\rm p} \bold{u}_{\rm p}})  =  -\phi_{\rm p} \nabla p + \nabla \cdot (\mu_{\rm p} \phi_{\rm p}(\nabla \bold{u}_{\rm p} + \nabla \bold{u}^\intercal_{\rm p})) \\
& \hspace{5.8em} + K(\bold{u_{\rm f}} - \bold{u_{\rm p}}) \nonumber, \\
\label{geep}
&\nabla \cdot ({\phi_{\rm p} \rho_{\rm p}  c_{\rm p}\bold{u}_{\rm p} T_{\rm p}})  =\nabla \cdot (\phi_{\rm p}\hat{k}_{\rm p}\nabla T_{\rm p}) + h(T_{\rm f} - T_{\rm p}).
\end{align}

The momentum exchange coefficient $K$ in Eqs.~(\ref{nsf}) and ~(\ref{nsp}) can be expressed as follows:
\begin{equation}
K = 
\begin{cases}
\begin{aligned}
&\frac{3}{4} C_{\rm d}  \frac{\phi_{\rm p} \phi_{\rm f} \rho_{\rm f} \lvert \bold{u}_{\rm f} - \bold{u}_{\rm p}\rvert} {d_{\rm p}} \phi_{\rm f}^{-2.65} , \text{ if $\phi_{\rm f} >$ 0.8} \\
&150\frac{\phi_{\rm p}^2 \mu_{\rm f}}{\phi_{\rm f}d_{\rm p}^2} + 1.75  \frac{\phi_{\rm p} \rho_{\rm f} \lvert \bold{u}_{\rm f} - \bold{u}_{\rm p}\rvert} {d_{\rm p}} , \text{ if $\phi_{\rm f} \leq $ 0.8},
\end{aligned}
\end{cases}
\end{equation}
where $d_{\rm p}$ is the particle diameter. As mentioned by Kolve~\cite{gidaspow1994multiphase}, the $C_{\rm d}$ is the drag coefficient, which is defined as follows:
\begin{equation}
\label{eqcd}
C_{\rm d} =
\begin{cases}
\begin{aligned}
&\dfrac{24}{Re_{\rm p}}  					,\quad \text{if $Re_{\rm p} < $ 1} 		\\
&\dfrac{24}{Re_{\rm p}}(1 + 0.15Re_{\rm p}^{0.687}) 	,\quad \text{if 1$ \leq Re_{\rm p} \leq $ 1000}.
\end{aligned}
\end{cases}
\end{equation}

The particle Reynolds number $Re_{\rm p}$ in Eq.~(\ref{eqcd}) is described as:
\begin{equation}
Re_{\rm p}  = \frac{\rho_{\rm f} \lvert \bold{u}_{\rm f} - \bold{u}_{\rm p}\rvert d_{\rm p}}{\mu_{\rm f}}.
\end{equation}

For monodispersed spherical particles, the heat exchange coefficient $h$ in Eqs.~(\ref{geef}) and~(\ref{geep}) can be expressed in the following equation.
\begin{equation}
h = \frac{6(1 - \phi_{\rm f})}{d_{\rm p}} h_{\rm p}.
\end{equation}

According to Wakao and Kaguei~\cite{wakao1982heat}, the fluid-particle heat transfer coefficient $h_{\rm p}$ can be obtained using the empirical correlation as follows:
\begin{equation}
Nu_{\rm p} = \frac{h_{\rm p} d_{\rm p}}{k_{\rm f}} = 2 + 1.1Re_{\rm p}^{0.6}Pr^{\frac{1}{3}},
\end{equation}
where $Pr$ is the Prandtl number for the fluid phase. Besides, as proposed by Kuipers et al.~\cite{kuipers1992numerical}, the effective thermal conductivities for both phases $\hat{k}_{\rm f}$ and $\hat{k}_{\rm p}$ are shown below:
\begin{align}
&\hat{k}_{\rm f} = \frac{\tilde{k}_{\rm f}}{\phi_{\rm f}},  \\
&\hat{k}_{\rm p} = \frac{\tilde{k}_{\rm p}}{\phi_{\rm p}}, 
\end{align}
where $\tilde{k}_{\rm f}$ and $\tilde{k}_{\rm p}$ can be further expressed as follows:
\begin{align}
&\tilde{k}_{\rm f} = \left(  1 - \sqrt{1 - \phi_{\rm f}} \right)k_{\rm f},  \\
\label{eqA}
&\tilde{k}_{\rm p} = \sqrt{1 - \phi_{\rm f}}\left(\omega A + (1-\omega)\Gamma \right)k_{\rm f}, \\
\label{eqB}
&\Gamma = \frac{2}{1 - \frac{B}{A}}\left\{ \frac{B(A-1)}{A(1-\frac{B}{A})^2} \ln\left(\frac{A}{B}\right) - \frac{B-1}{(1-\frac{B}{A})} - \frac{B+1}{2} \right\},
\end{align}
where $A$ and $B$ in Eqs.~(\ref{eqA}) and~(\ref{eqB}) are given by the equations $A = \frac{k_{\rm p}}{k_{\rm f}}$ and $B = 1.25\left( \frac{1 - \phi_{\rm f}}{\phi_{\rm f}} \right)^{\frac{10}{9}}$, respectively. According to Kalteh et al.~\cite{kalteh2011eulerian}, the value of $\omega$ is chosen as $7.26 \times 10^{-3}$ in this study. 

Fig.~\ref{Fig.1} shows the analysis domain and dimensions for the optimization problems.  To account for the effect of pumping power on heat transfer performance, we fixed the pressure drop between the inlet and the outlet in the numerical examples. The boundary conditions for state variables in the first numerical example are shown below:
\begin{alignat}{4}
&p = p^{\rm in}   																&\text{\quad  on \quad}&\mathit{\Gamma}^{\rm in}, 		\\
&T_{\rm f} = T^{\rm in}_{\rm f}, \quad T_{\rm p} = T^{\rm in}_{\rm p} 			&\text{\quad  on \quad}&\mathit{\Gamma}^{\rm in}, 		\\
&p = p^{\rm out}  																&\text{\quad  on \quad}&\mathit{\Gamma}^{\rm out},		\\
&\bold{\phi}_{\rm p} = \phi^{\rm in}_{\rm p}, 												&\text{\quad  on \quad}&\mathit{\Gamma}^{\rm in}, 		\\
&\bold{u}_{\rm f} = \bold{0}, \quad \bold{u}_{\rm p} = \bold{0}			 						&\text{\quad  on \quad}&\mathit{\Gamma}^{\rm w}, 		\\
&\nabla {\phi}_{\rm p} \cdot \bold{n} = 0												&\text{\quad  on \quad}&\mathit{\Gamma}^{\rm w}, \mathit{\Gamma}^{\rm out}, \\
&\nabla {T}_{\rm f} \cdot \bold{n} = 0, 	\quad \nabla {T}_{\rm p} \cdot \bold{n} = 0					&\text{\quad  on \quad}&\mathit{\Gamma}^{\rm w}, \mathit{\Gamma}^{\rm out}.
\end{alignat}

In addition, we apply the symmetric boundary condition in the first numerical example to reduce computational cost. In the second numerical example, since we optimize the flow field in the unit cell, free-slip boundary conditions are applied to the walls. The types of the remaining boundary conditions are the same as those in the first numerical example.
\section{Problem formulation}
\label{Problem formulation}
\subsection{Representation of solid and fluid regions}
\label{Representation of solid and fluid regions}
Topology optimization aims to enhance the objective function by optimizing the material distribution within the design domain $D$. Initially, we define a characteristic function for the representation of solid and fluid regions to determine the material distribution, as follows:
\begin{equation}
\bold{\chi}(\bold{x}) = 
\begin{cases}
1 , \text{ if $\bold{x} \in \mathit{\Omega}$}\\
0, \text{ if $\bold{x} \in D \backslash \mathit{\Omega}$},
\end{cases}
\end{equation}
where $\bold{x}$ is the position in $D$, and $\mathit{\Omega}$ represents the flow domain in $D$. The characteristic function $\bold{\chi}(\bold{x})$ is inherently discontinuous, which necessitates the relaxation techniques to facilitate a smooth transition for numerical analysis. In this study, we employ a continuous function, $\gamma(\bold{x})$, commonly used in the density method, with values ranging $0 \leq \gamma(\bold{x}) \leq 1$.

Checkerboard patterns may arise in the material distribution due to numerical instabilities~\cite{diaz1995checkerboard}. To address this issue, we implement convolution filtering as outlined below~\cite{bourdin2001filters}:
\begin{equation}
\label{filter}
\tilde{\gamma}_{k} = \frac{\sum_{i \in N_{e, k}}w(\bold{x}_{\mathit{i}})\gamma_{\mathit{i}}}{\sum_{i \in N_{e, k}}w(\bold{x}_{\mathit{i}})}.
\end{equation}

Here, $\tilde{\gamma}_{k}$ represents the filtered design variable on node $k$, $N_{e, k}$ denotes the set of nodes located within the filtering radius $R$, and $ w(\bold{x}_{\mathit{i}})$ is the weighting function defined as below.
\begin{equation}
w(\bold{x}_{\mathit{i}}) = \mathit{R} - \lvert \bold{x}_{\mathit{i}} - \bold{x}_{\mathit{k}}\rvert.
\end{equation}

Although filtering can address the checkerboard pattern issue, it often results in the presence of gray-scale design variables in the material distribution. To address this issue, we employ the projection function to obtain the optimized configuration closer to a binary distribution as follows~\cite{wang2011projection}:
\begin{equation}
\label{projection}
\bar{\tilde{\gamma}}_{k} = \frac{\tanh(\beta\eta) + \tanh(\beta(\tilde{\gamma}_{k} - \eta))}{\tanh(\beta\eta) + \tanh(\beta(1 - \eta))},
\end{equation}
where $\beta$ is the projection parameter, $\eta$ is the threshold parameter and $\bar{\tilde{\gamma}}_{k}$ is the design variable after projection.

\subsection{Governing equations with fluid and solid representations}
In Section~\ref{Representation of solid and fluid regions}, we defined the design variables used to represent the solid and fluid regions. To determine the distribution of design variables during the optimization process, we further introduced fictitious forces and heat sources dependent on the design variables and incorporated them into the equations. The re-formulated equations are presented as follows:
\begin{align}
\label{rensf}
&\nabla \cdot ({\phi_{\rm f} \rho_{\rm f} \bold{u}_{\rm f} \bold{u}_{\rm f}}) &&\negthickspace\negthickspace\negthickspace =  -\phi_{\rm f} \nabla p + \nabla \cdot (\mu_{\rm f} \phi_{\rm f}(\nabla \bold{u}_{\rm f} + \nabla \bold{u}^\intercal_{\rm f})) \\
&  &&\negthickspace\negthickspace\negthickspace + K(\bold{u_{\rm p}} - \bold{u_{\rm f}}) - \phi_{\rm f}\alpha \bold{u}_{\rm f} \nonumber, \\[10pt]
\label{regeef}
&\nabla \cdot ({\phi_{\rm f} \rho_{\rm f}  c_{\rm f}\bold{u}_{\rm f} T_{\rm f}}) &&\negthickspace\negthickspace\negthickspace =\nabla \cdot (\phi_{\rm f}k_{\rm eff, f}\nabla T_{\rm f}) + h(T_{\rm p} - T_{\rm f}) \\
&  &&\negthickspace\negthickspace\negthickspace + \phi_{\rm f}Q (T_{\rm Q} - T_{\rm f}) \nonumber, \\[10pt]
\label{rensp}
&\nabla \cdot ({\phi_{\rm p} \rho_{\rm p} \bold{u}_{\rm p} \bold{u}_{\rm p}}) &&\negthickspace\negthickspace\negthickspace =  -\phi_{\rm p} \nabla p + \nabla \cdot (\mu_{\rm p} \phi_{\rm p}(\nabla \bold{u}_{\rm p} + \nabla \bold{u}^\intercal_{\rm p})) \\
&  &&\negthickspace\negthickspace\negthickspace + K(\bold{u_{\rm f}} - \bold{u_{\rm p}}) - \phi_{\rm p}\alpha \bold{u}_{\rm p} \nonumber, \\[10pt]
\label{regeep}
&\nabla \cdot ({\phi_{\rm p} \rho_{\rm p}  c_{\rm p}\bold{u}_{\rm p} T_{\rm p}}) &&\negthickspace\negthickspace\negthickspace =\nabla \cdot (\phi_{\rm p}k_{\rm eff, p}\nabla T_{\rm p}) + h(T_{\rm f} - T_{\rm p}) \\
&  &&\negthickspace\negthickspace\negthickspace + \phi_{\rm p}Q (T_{\rm Q} - T_{\rm p}) \nonumber.
\end{align}

The inverse permeability $\alpha$ and volumetric heat generation coefficient $Q$ can be described as follows~\cite{matsumori2013topology}:
\begin{align}
&\alpha(\bar{\tilde{\gamma}}) = \overline{\alpha} + (\underline{\alpha} - \overline{\alpha})\bar{\tilde{\gamma}} \frac{1 + q}{\bar{\tilde{\gamma}} + q}, \\
&Q(\bar{\tilde{\gamma}}) = \overline{Q}(1 - \bar{\tilde{\gamma}}),
\end{align}
where $\overline{\alpha}$ and $\underline{\alpha}$ are the maximum and minimum values of the inverse permeability, respectively, $\overline{Q}$ is the maximum value of the volumetric heat generation coefficient, and $q$ is the convex parameter. It is worth noting that $\overline{Q}$ is also necessary when converting a three-dimensional problem into a two-dimensional problem. By employing methods such as the two-layer model~\cite{yan2019topology}, results comparable to those obtained in three dimensions can be achieved if the height of the flow path is sufficiently small. However, for simplicity, we follow previous research~\cite{matsumori2013topology, yaji2016topology} and determine $\overline{Q}$ based on prior parameter investigations. By incorporating fictitious forces and heat sources into the equations, we can evaluate the influence of design variables on the objective function through sensitivity analysis.

\subsection{Heat transfer maximization problem}
Heat transfer between the solid and fluid regions is a critical factor in the design of MCHS. In this optimization problem, since the pressure drop between the inlet and outlet is fixed, the objective function only needs to account for the heat transfer between the solid and the fluid. Therefore, the optimization problem is formulated as a heat generation maximization problem as follows:
\begin{equation}
\begin{split}
&\text{maximize }   J = \int_{\mathit{D}} Q(T_{\rm Q} - T_{\rm ave})d\mathit{\Omega}, \\
&0 \leq \gamma(\bold{x}) \leq 1 \text{ for } \forall \bold{x} \in D,
\end{split}
\end{equation}
where $T_{\rm Q}$ represents the temperature in the solid region, $T_{\rm ave}$ represents the average temperature of the fluid phase and the particle phase, expressed as follows~\cite{boulet2002influence}, where $l$ denotes the phase with $l \in \{ \rm{f},\rm{p}\}$.
\begin{equation}
T_{\rm ave} = \frac{\sum\limits_{l \in \{ \rm f, p \}}\left( \int_{\mathit{D}} \rho_{l} \left|\bold{u}_{l}\right| c_{l}T_{l}d\mathit{\Omega} \right)}{\sum\limits_{l \in \{ \rm f, p \}}\left( \int_{\mathit{D}} \rho_{l} \left|\bold{u}_{l}\right| c_{l}d\mathit{\Omega}  \right)}
\end{equation}
\section{Numerical implementation}
\label{Numerical implementation}
We develop a custom finite difference flow solver incorporating the NAPPLE (Nonstaggered Artificial Pressure for Pressure-Linked Equation) algorithm. This algorithm allows state variables to be defined on a collocated grid under specific assumptions. Based on these assumptions, we can derive the pressure-linked equation and discretize it using the harmonic scheme to eliminate the checkerboard pressure pattern. For further details, refer to~\cite{chen2024topology, lee1992artificial}. Besides, the momentum and energy equations are discretized for both phases using a weighting function scheme~\cite{shong1989weighting}. In fluid simulation, we use the maximum residual values among all grid points as the convergence criteria expressed by the following equations:
\begin{align}
&\mathop{\max} \left| u^{i} - u^{i - 1} \right| \leq \rm{tol}_\mathit{u}, \\[10pt]
&\mathop{\max} \left| v^{i} - v^{i - 1} \right| \leq \rm{tol}_\mathit{v}, \\[10pt]
&\mathop{\max} \left| T^{i} - T^{i - 1} \right| \leq \rm{tol}_\mathit{T},
\end{align}
where $i$ represents iteration count, $u$ and $v$ are the velocity components in the $x$- and $y$-directions, respectively. The convergence criteria for the state variables $u$, $v$ and $T$ for both phases are denoted as $\rm{tol}_\mathit{u}$, $\rm{tol}_\mathit{v}$ and $\rm{tol}_\mathit{T}$ and are set to $1 \times 10^{-6}$, respectively.

\begin{figure}[h]
\centering
\includegraphics[width=0.5\textwidth]{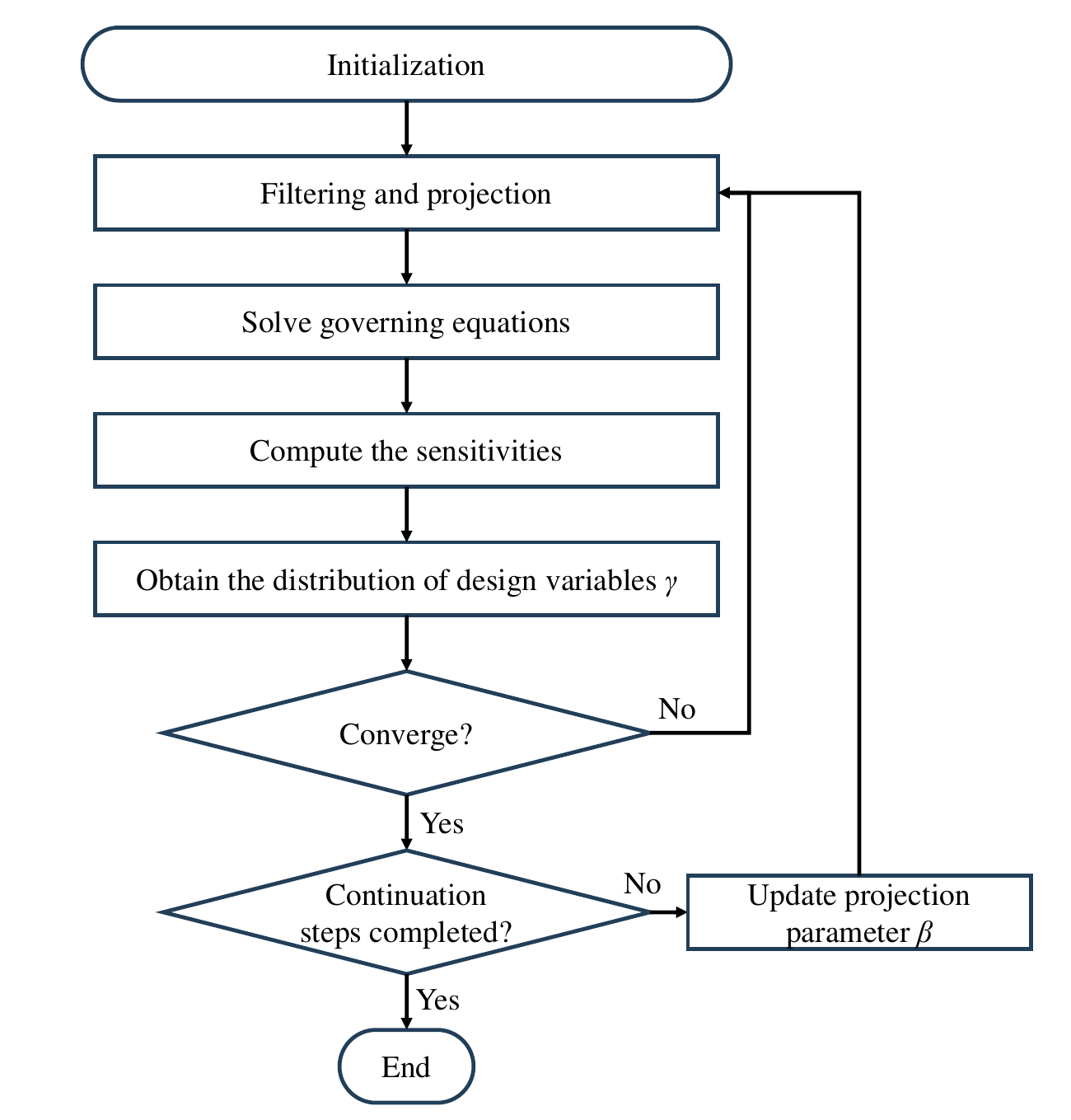} 
\caption{Flowchart of the topology optimization process.}
\label{Fig.2}
\end{figure}

Fig.~\ref{Fig.2} shows the flowchart of the optimization process. In the first step, we initialize all the state and design variables. In the second step, we apply the filter in Eq.~(\ref{filter}) and projection in Eq.~(\ref{projection}) on the design variables. In the third and fourth steps, the numerical simulation is performed, and the sensitivities are obtained using the automatic differentiation with Adept C++~\cite{hogan2014fast}. In the fifth step, we use the GCMMA algorithm~\cite{svanberg2002class} to determine the distribution of design variables. In the optimization process, we adopt a continuation method, gradually increasing the value of $\beta$ in Eq.~(\ref{projection}) to remove gray-scale design variables. If the optimization meets the convergence criterion, the value of $\beta$ is doubled, and the process returns to the second step. This loop continues until all continuation steps are completed. If the criterion is unmet, the process returns to the second step without modifying $\beta$. The convergence criterion for the optimization is as follows:
\begin{equation}
\left| \frac{ J_{n} - J_{n-1}}{J_{n}} \right| < \epsilon.
\end{equation}

Here, $n$ represents the number of optimization iterations, $\epsilon$ represents the criterion for the optimization convergence and is set to $1 \times 10^{-4}$ in this study.
\begin{figure}[h]
\centering
\includegraphics[width=0.45\textwidth]{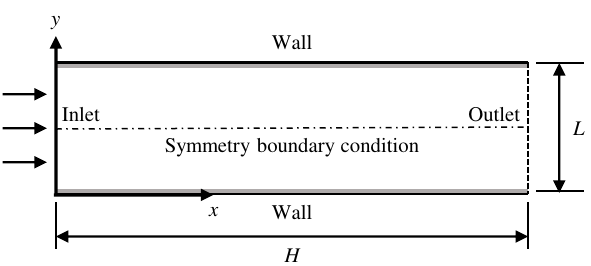} 
\caption{Analysis domain for the numerical scheme validation.}
\label{Fig.3}
\end{figure}
\section{Results and discussion}
\label{Results and discussion}
\subsection{Validation of numerical method}
We validate the feasibility of the numerical scheme through a numerical example from prior research~\cite{kalteh2011eulerian}. In this problem, the nanofluid flows between two parallel plates separated by $L$ = 200 $\mu$m. Each plate has a length 100 times its inlet width ($H/L = 100$), as shown in Fig.~\ref{Fig.3}. The analysis domain is discretized with a quadrilateral mesh of $2000 \times 30$. The nanofluid consists of water as the base fluid and copper as the nanoparticle material. At the inlet, the nanofluid enters with a uniform velocity and a temperature of 293 K, while the static pressure at the outlet is set to zero. The side walls are subjected to no-slip boundary conditions, with the wall temperature set to 303 K. In addition, the symmetry boundary condition is utilized for the numerical analysis. The material properties of the nanofluid are summarized in the table~\ref{validation1}.

\begin{figure}[h]
\centering
\includegraphics[width=0.5\textwidth]{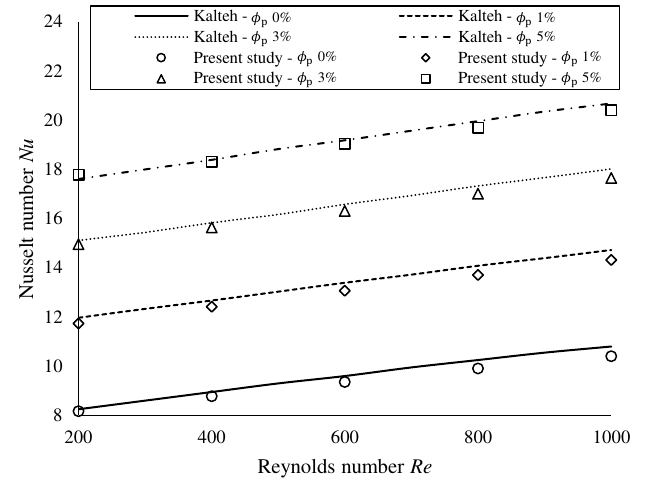} 
\caption{Average Nusselt number is evaluated with the Reynolds number at different particle volume fractions in the analysis domain shown in Fig.~\ref{Fig.3}.}
\label{Fig.4}
\end{figure}

In this study, we first compare the average Nusselt number variations with the Reynolds number at different particle volume fractions. The definitions and calculations of Reynolds and Nusselt numbers are based on those provided by Kalteh et al.~\cite{kalteh2011eulerian}. As illustrated in Fig.~\ref{Fig.4}, the numerical results align well with those reported in the previous study. An increase in Reynolds number enhances the convective heat transfer effect, while a higher particle volume fraction improves the effective thermal conductivity of the nanofluid. Consequently, the average Nusselt number increases accordingly.

Additionally, we examined the average Nusselt number for different particle diameters under $Re$ = 500 and $\phi_{\rm p}$ = 1$\%$. As reported in Kalteh et al.~\cite{kalteh2011eulerian}, the ratio of the Nusselt number for nanofluids to water is 40.43$\%$ and 42.47$\%$ for $d_{\rm p}$ = 100 nm and $d_{\rm p}$ = 30 nm, respectively, with the Nusselt number of water reported as 9.3. Based on these values, the average Nusselt numbers for $d_{\rm p}$ = 100 nm and $d_{\rm p}$ = 30 nm can be calculated as 13.06 and 13.25, respectively. Under the same conditions, the average Nusselt numbers for $d_{\rm p}$ = 100 nm and $d_{\rm p}$ = 30 nm in this study are 12.74 and 13.42, showing a deviation of less than 2.5$\%$ compared to the prior study.

\begin{table}[h]\centering
\caption{Material properties of the Cu-water nanofluid.}
\label{validation1}
\fontsize{4}{4}\selectfont
\resizebox{1.0\linewidth}{!}{
\begin{tabular}{@{}lccll@{}}
																												 \specialrule{.1em}{.05em}{.05em} 
Parameters			& \multicolumn{1}{l}{Symbols} 		& \multicolumn{2}{c}{Values}   		& Unit 	\rule[1.8ex]{0pt}{1.8ex}\\ 	\specialrule{.1em}{.05em}{.05em} 
Fluid viscosity                  & $\mu_{\rm f}$     	      		& \multicolumn{2}{c}{0.001}      	& $\rm Pa\cdot s$  			\rule[1.8ex]{0pt}{1.8ex} \\
Fluid density                  	& $\rho_{\rm f}$     	      		& \multicolumn{2}{c}{1000}      	& $\rm kg/m^{3}$  		\rule[1.8ex]{0pt}{1.8ex} \\
Fluid heat capacity           & $c_{\rm f}$     	      			& \multicolumn{2}{c}{4181.8}      	& $\rm J / kg \, K$  		\rule[1.8ex]{0pt}{1.8ex} \\
Fluid thermal conductivity        	& $k_{\rm f}$     	      		& \multicolumn{2}{c}{0.6}      		& $\rm W / m \, K$  		\rule[1.8ex]{0pt}{1.8ex} \\
Particle viscosity              & $\mu_{\rm p}$     	      		& \multicolumn{2}{c}{$1.38 \times 10^{-3}$}      	& $\rm Pa\cdot s$  		\rule[1.8ex]{0pt}{1.8ex} \\
Particle density             	& $\rho_{\rm p}$  				& \multicolumn{2}{c}{8954}   				& $\rm kg/m^{3}$  	\rule[1.8ex]{0pt}{1.8ex}\\
Particle heat capacity       & $c_{\rm p}$     	      			& \multicolumn{2}{c}{383.1}      			& $\rm J / kg \, K$  	\rule[1.8ex]{0pt}{1.8ex} \\
Particle thermal conductivity        	& $k_{\rm p}$     	      	& \multicolumn{2}{c}{386}      				& $\rm W / m \, K$  			\rule[1.8ex]{0pt}{1.8ex} \\
Particle diameter           		& $d_{\rm p}$       		& \multicolumn{2}{c}{100}    				& nm    		 			\rule[1.8ex]{0pt}{1.8ex}\\ \specialrule{.1em}{.05em}{.05em} 
\end{tabular}}
\end{table}
\subsection{Heat transfer optimization problem}
\subsubsection{Details of numerical setting}
\label{Details of numerical setting}
\begin{figure}[h]
\centering
\includegraphics[width=0.45\textwidth]{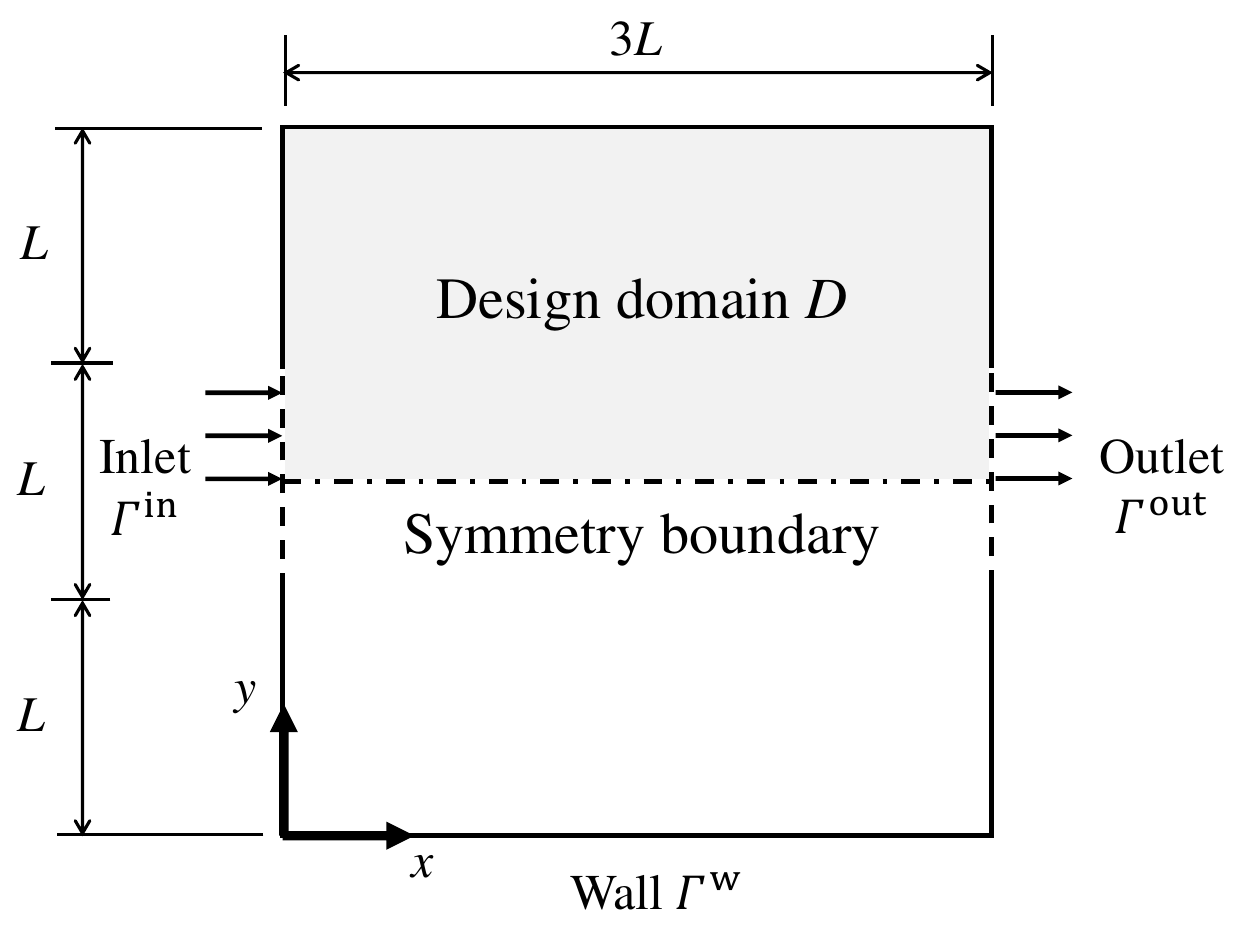} 
\caption{Dimensions of the design domain $D$ utilized for the heat transfer optimization problem.}
\label{Fig.15}
\end{figure}

In the first numerical example, the default material properties for the nanofluid are listed in Table~\ref{validation1}. Fig.~\ref{Fig.15} shows the schematic diagram of the analysis domain in this problem, with the inlet length $L$ set to 30 $\mu$m. Although the overall size of the analysis domain is relatively small for MCHS, this configuration is commonly used in topology optimization problems for heat sinks~\cite{matsumori2013topology, yaji2016topology}. As a preliminary validation case, this setup allows for evaluating the influence of different parameters on the optimized flow field at a reasonable computational cost and verifies whether the proposed method can generate flow fields with superior performance.

In the fluid simulation setup, the analysis domain is discretized using a quadrilateral mesh of  $100 \times 50$. The pressure drop is set to 5000 $ \rm Pa$, and the inlet particle volume fraction $\phi^{\rm in}_{\rm p}$ is $1 \%$. The corresponding fluid Reynolds number is 67.1. For the temperature boundary conditions, the inlet temperature $T^{\rm in}$ for both phases is set to 300 K, while solid temperature $T_{\rm Q}$ is set to 360 K. 

In the optimization setup, the initial design variable is set to 0.9, while the inverse permeability $\alpha$ and the volumetric heat generation coefficient $Q$ have maximum values of $4.44 \times 10^{10}$ and $5.33 \times 10^{10}$, respectively. The convex parameter $q$ is set to 0.05. In the filtering and projection, the filter radius $R$ is set to 1.8 $\rm \mu m$, and the maximum value of the projection parameter $\beta_{\rm max}$ is set to 64.
\begin{figure*}[htbp]
\centering
\includegraphics[width=0.8\textwidth]{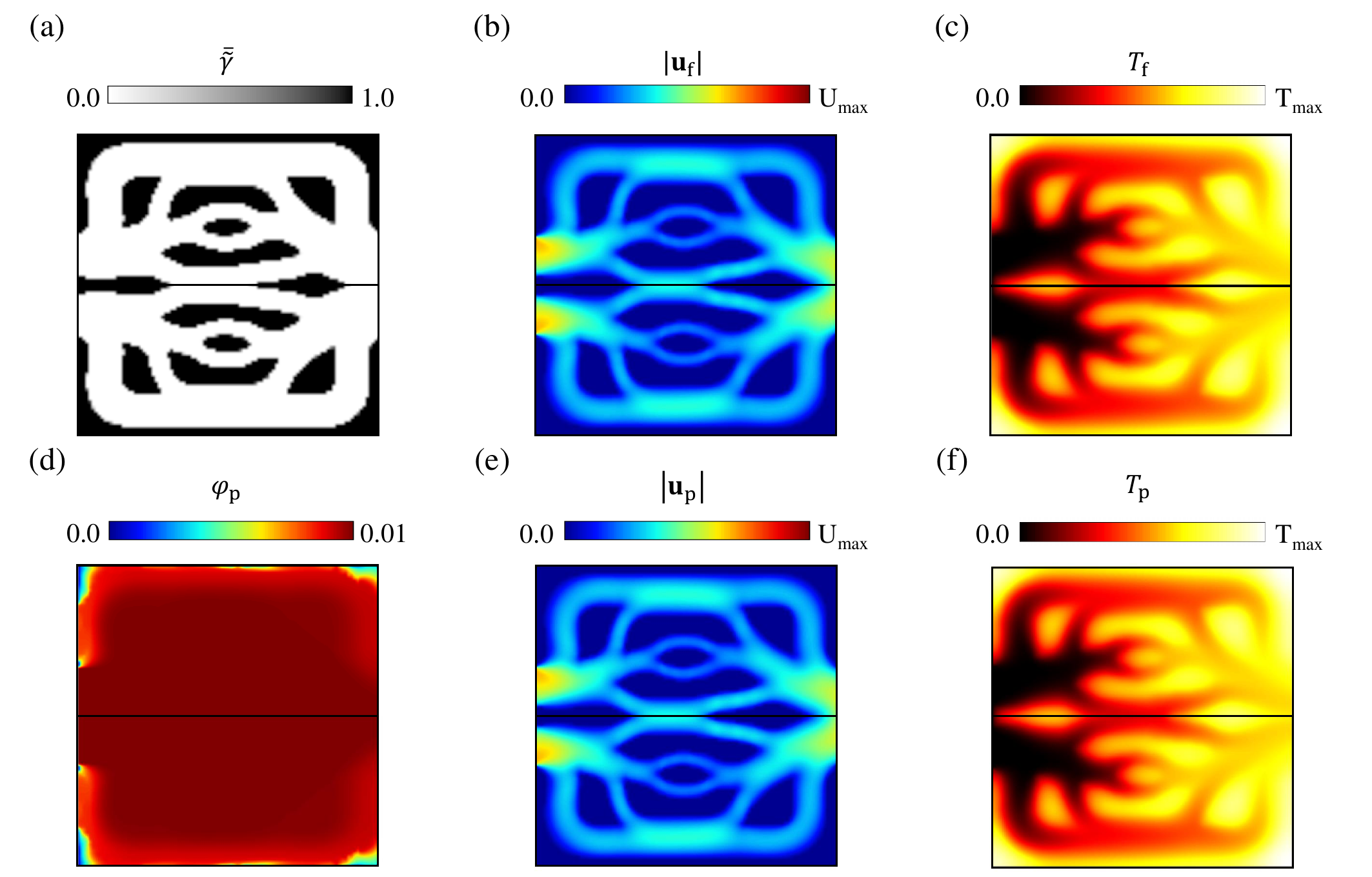} 
\caption{Optimized Results and distribution of physical quantities with: (a) optimized flow field (b) fluid velocity magnitude(c) fluid phase temperature (d) particle volume fraction (e) particle velocity magnitude (f) particle phase temperature. The maximum values on the color bar for velocity and temperature for both phases are 1.77 $\rm m/s$ and 355.8 $\rm K$, respectively.}
\label{Fig.5}
\end{figure*}
\subsubsection{Optimized flow field}
The optimized flow field is obtained using the default parameters outlined in Section~\ref{Details of numerical setting}. Fig.~\ref{Fig.5} illustrates the optimized flow field alongside the velocity and temperature distributions for both phases. From Fig.~\ref{Fig.5}(c) and (f), it can be observed that as the nanofluid flows around the heat sink fins, it absorbs heat generated by the solid through the thermal exchange, thereby achieving a cooling effect. Due to the interactions between the fluid and particles, the velocity and temperature distributions of the fluid phase and particle phase are nearly identical. Notably, as shown in Fig.~\ref{Fig.5}(d), the particle volume fraction is distributed uniformly across the design domain, even within the solid regions. This is due to the limitations of the fictitious body force formulation, which causes particles to remain inside the solid regions. However, since the particle velocity within the solid regions is extremely low, its influence on the optimization results is minimal. 

\begin{figure}[h]
\centering
\includegraphics[width=0.4\textwidth]{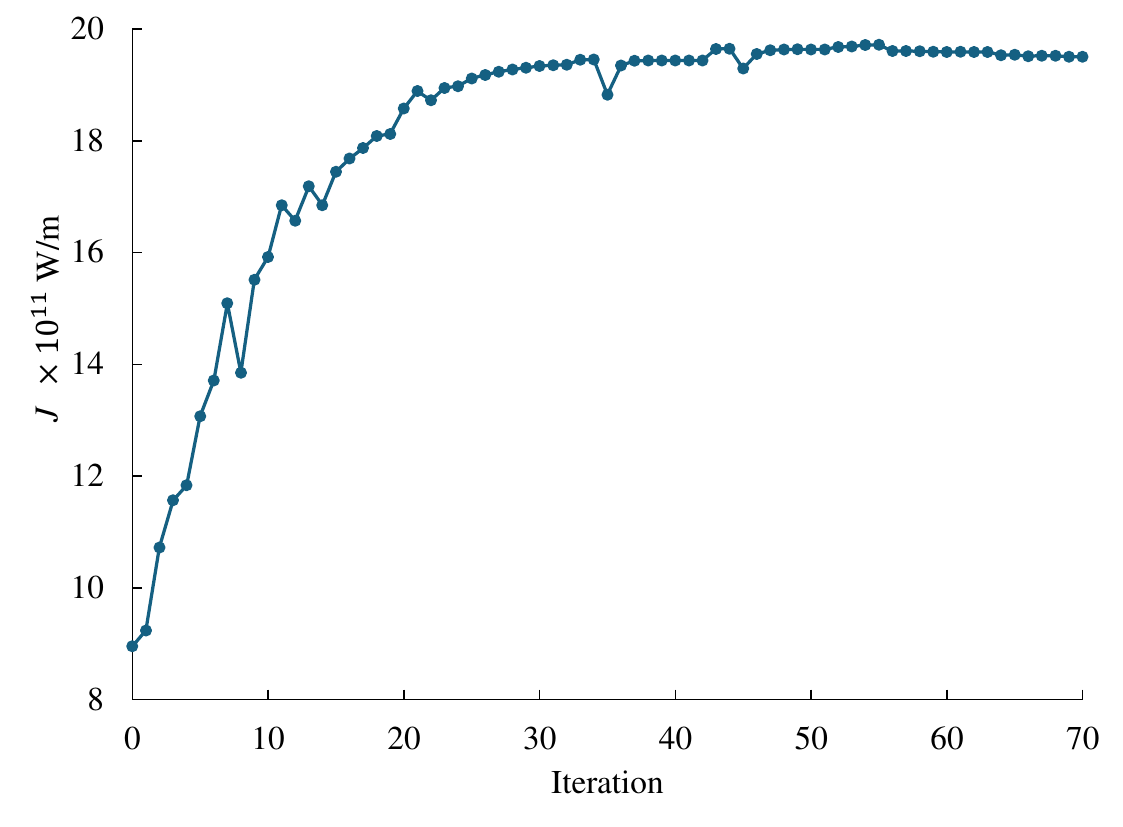} 
\caption{Convergence history of the objective function $J$ for the optimized flow field shown in Fig.~\ref{Fig.5}.}
\label{Fig.6}
\end{figure}

Fig.~\ref{Fig.6} illustrates the convergence history for the optimized flow field, showing a progressive increase in the objective function value with each iteration. However, the objective function value may occasionally experience a sudden drop. This is because the projection function in Eq.~(\ref{projection}) imposes abrupt changes to the design variable field, resulting in short-term performance reduction.
\begin{figure*}[htbp]
\centering
\includegraphics[width=0.8\textwidth]{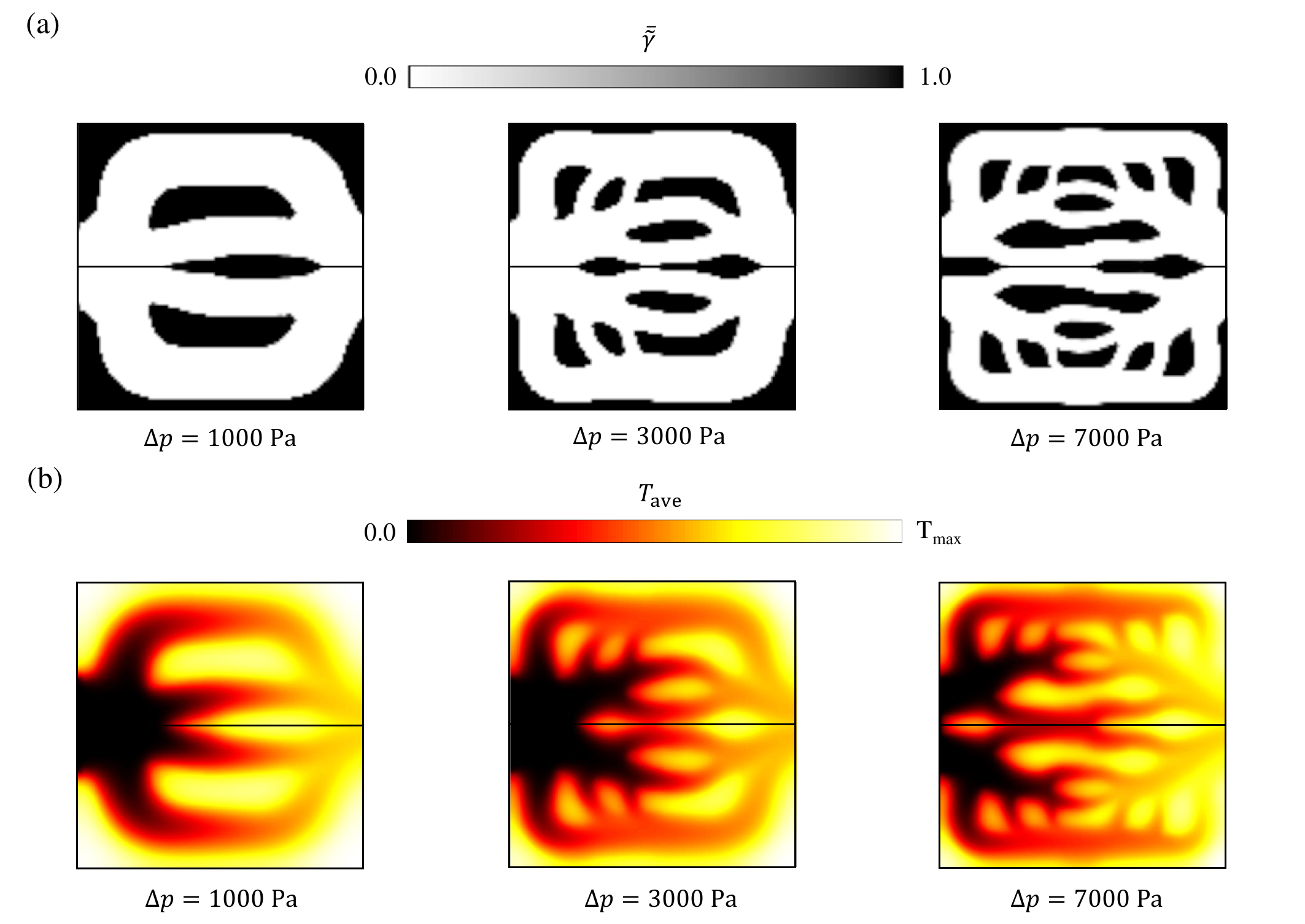} 
\caption{Optimized Results and distribution of physical quantities with different pressure drop: (a) optimized flow field (b) average temperature. The maximum values on the color bar for the average temperature are 359.3 $\rm K$,  358.3 $\rm K$, and 354.3 $\rm K$ with $\Delta p$ = 1000 Pa, 3000 Pa and 7000 Pa, respectively.}
\label{Fig.7}
\end{figure*}
\subsubsection{Effect of pressure drop}
First, we examine the effects of different pressure drops on the optimized configurations, as shown in Fig.~\ref{Fig.7}. The results indicate that the number of branches in the optimized flow field increases as the pressure drop increases. An increase in pressure drop enhances the overall convective effects within the flow field, promoting the formation of additional flow branches. Although increased branching leads to greater flow resistance, it also enhances heat dissipation performance by providing additional heat transfer surfaces and disrupting the thermal boundary layer. Table~\ref{crosscheck1} presents the results of the crosscheck. It can be observed that the flow fields optimized under specific pressure drop conditions exhibit the best performance at their corresponding pressure drops. This demonstrates that the proposed method can generate optimized flow fields under different pressure drop conditions.
\setlength{\arrayrulewidth}{0.1mm}
\begin{table}[h]
\caption{Crosscheck of objective function $ J$ ($\times 10^{11}$) with pressure drop.}
\label{crosscheck1}
\renewcommand{\arraystretch}{1.5}
\fontsize{4}{4}\selectfont
\resizebox{1.0\linewidth}{!}{
\begin{tabular}{lllll}
\hline
\multirow{2}{*}{Analysis $\Delta p$} & \multicolumn{4}{c}{Optimization $\Delta p$} \\ \cline{2-5} 
                          	& 1000   & 3000   & 5000   & 7000  \\ \cline{1-1}
1000                      & 11.0   & 9.3     & 7.4     & 6.3  \\
3000                      & 14.6   & 16.4   & 15.6   & 14.6  \\
5000                      & 15.8   & 19.0   & 19.3   & 19.1  \\
7000                      & 16.4   & 20.4   & 21.8   & 21.9  \\ \hline
\end{tabular}
}
\end{table}
\begin{figure*}[htbp]
\centering
\includegraphics[width=0.8\textwidth]{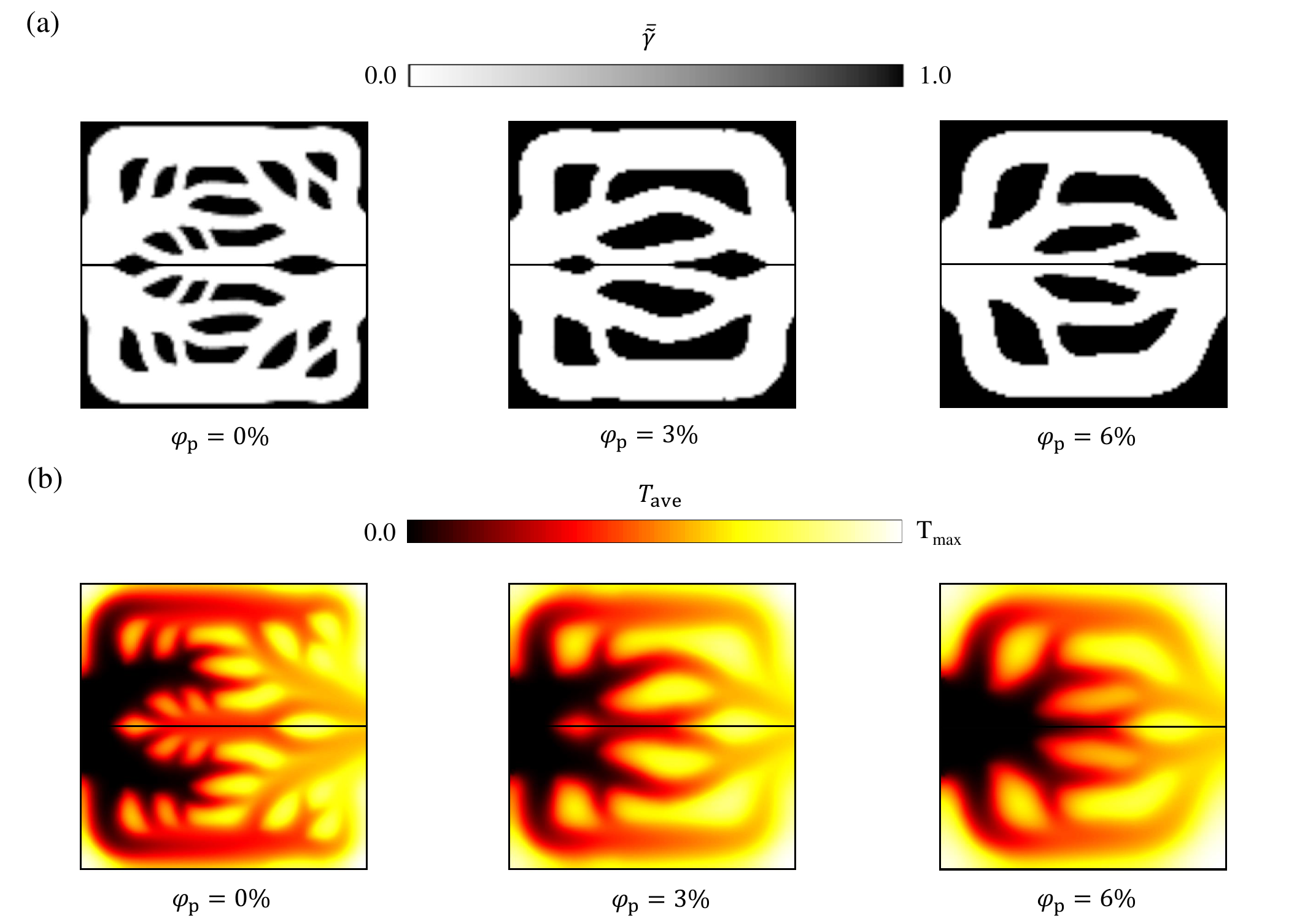} 
\caption{Optimized Results and distribution of physical quantities with different particle volume fraction: (a) optimized flow field (b) average temperature. The maximum values on the color bar for the average temperature are 354.6 $\rm K$,  356.0 $\rm K$, and 357.9 $\rm K$ with $\phi_{\rm p}$ = 0 $\%$, 3$\%$ and 6$\%$, respectively.}
\label{Fig.8}
\end{figure*}
\begin{figure*}[htbp]
\centering
\includegraphics[width=0.8\textwidth]{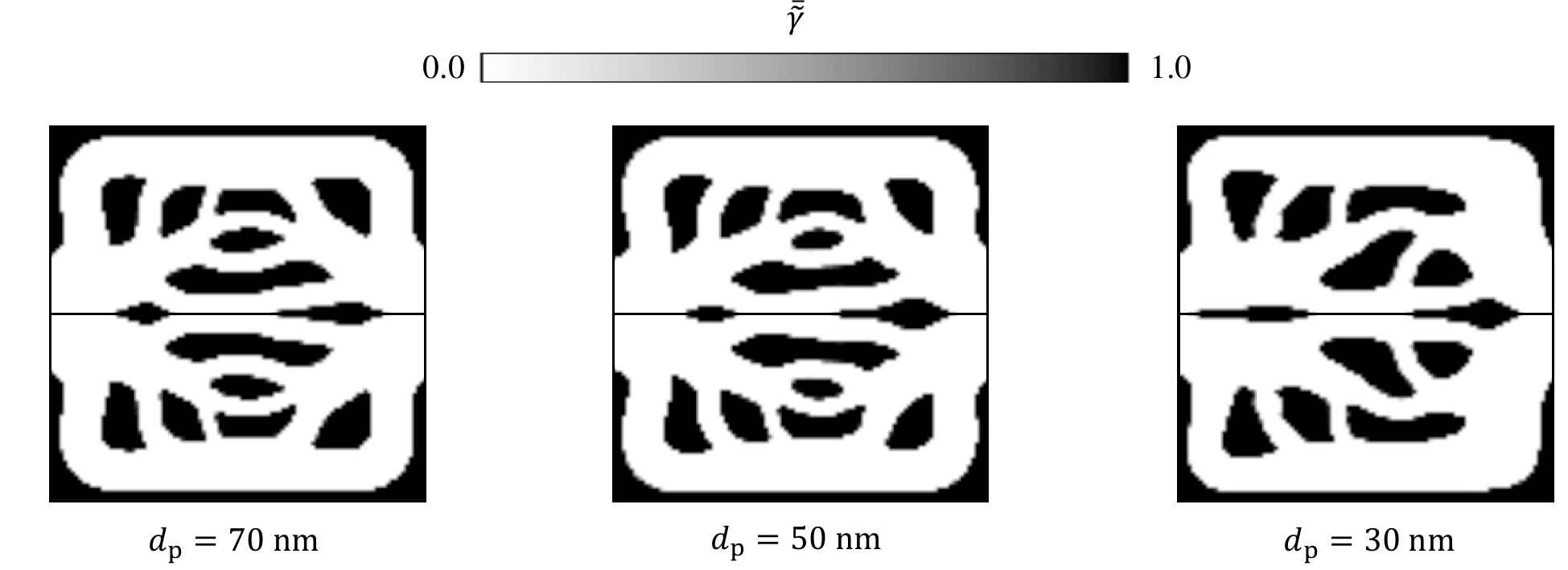} 
\caption{Optimized flow fields with different particle diameters.}
\label{Fig.13}
\end{figure*}
\subsubsection{Effect of particle volume fraction}
Fig.~\ref{Fig.8} presents the optimized configurations for various particle volume fractions. As the particle volume fraction increases, the number of branches in the optimized flow field decreases. An increase in the particle volume fraction also raises the fluid viscosity, leading to greater flow resistance. Consequently, reducing the number of branches minimizes flow losses under a fixed pressure drop. This adjustment allows nanofluids with higher particle volume fractions to sustain sufficient flow velocity, thereby improving convective heat transfer performance. Table~\ref{crosscheck2} presents the crosscheck results at different particle volume fractions. Under a fixed pressure drop, the flow field optimized at $\phi_{\rm p}$ = 6$\%$ exhibits the best performance. Although the performance of each optimized flow field generally improves with higher particle volume fractions, the flow field optimized for a specific fraction consistently outperforms those optimized for other fractions.
\setlength{\arrayrulewidth}{0.1mm}
\begin{table}[h]
\caption{Crosscheck of objective function $J$ ($\times 10^{11}$) with volume fraction.}
\label{crosscheck2}
\renewcommand{\arraystretch}{1.5}
\fontsize{4}{4}\selectfont
\resizebox{1.0\linewidth}{!}{
\begin{tabular}{lllll}
\hline
\multirow{2}{*}{Analysis $\phi_{\rm p}$} & \multicolumn{4}{c}{Optimization $\phi_{\rm p}$} \\ \cline{2-5} 
                          	& 0$\%$   & 1$\%$   & 3$\%$   & 6$\%$  \\ \cline{1-1}
0$\%$                      & 18.2   & 18.1   & 17.5   & 17.0  \\
1$\%$                      & 19.2   & 19.6   & 19.4   & 19.1  \\
3$\%$                      & 19.5   & 20.2   & 20.4   & 20.2  \\
6$\%$                      & 19.8   & 20.5   & 21.0   & 21.1  \\ \hline
\end{tabular}
}
\end{table}
\subsubsection{Effect of particle diameter}
Fig.~\ref{Fig.13} illustrates the optimized flow fields for different particle diameters. The topological structures of the optimized flow fields remain unchanged for $d_{\rm p}$ = 70 nm and 50 nm, with only minor variations in shape. Similarly, the optimized flow fields for $d_{\rm p}$ = 100 nm and 30 nm show slight differences in their topological structures. In summary, particle diameter has a minor effect compared to the influence of particle volume fraction on the optimized flow field. As reported by Zhang et al.~\cite{zhang2020topology}, a single-phase model was used to examine the effect of particle diameter, and it similarly showed that particle diameter has a relatively minor impact on the optimized flow field.

In addition, as shown in the crosscheck results in Table~\ref{crosscheck3}, the differences in the objective function values among the particle diameters are minimal for each optimized flow field. This is because, within the conditions explored in this study, the relatively low Reynolds number and particle volume fraction minimize the influence of particle diameter variation on nanofluid behavior. In the prior work, Kalteh et al.~\cite{kalteh2011eulerian} employed an Eulerian-Eulerian model to investigate the effect of particle diameter on the Nusselt number enhancement at $Re$ = 500 and $\phi_{\rm p}$ = 1$\%$, using the same nanofluid and the analysis domain shown in Fig.~\ref{Fig.3}. Their results indicated that even at $Re$ = 500, reducing the particle diameter from 100 nm to 30 nm yielded only a 2$\%$ improvement. Given the lower Reynolds number with the same particle volume fraction in this study, the effect of particle diameter is further diminished. Consequently, the configurations of the optimized flow fields exhibit little variation, and their performance across different particle diameters is nearly identical.
\setlength{\arrayrulewidth}{0.1mm}
\begin{table}[h]
\caption{Crosscheck of objective function $ J $($\times 10^{11}$) with particle diameter.}
\label{crosscheck3}
\renewcommand{\arraystretch}{1.5}
\fontsize{4}{4}\selectfont
\resizebox{1.0\linewidth}{!}{
\begin{tabular}{lllll}
\hline
\multirow{2}{*}{Analysis $d_{\rm p}$} & \multicolumn{4}{c}{Optimization $d_{\rm p}$} \\ \cline{2-5} 
                          & 30   & 50   & 70   & 100  \\ \cline{1-1}
30                      & 18.9   & 19.2   & 19.2   & 19.3  \\
50                      & 19.0   & 19.4   & 19.4   & 19.5  \\
70                      & 19.0   & 19.4   & 19.4   & 19.5  \\
100                    & 19.0   & 19.3   & 19.4   & 19.6  \\ \hline
\end{tabular}
}
\end{table}
\begin{figure}[h]
\centering
\includegraphics[width=0.45\textwidth]{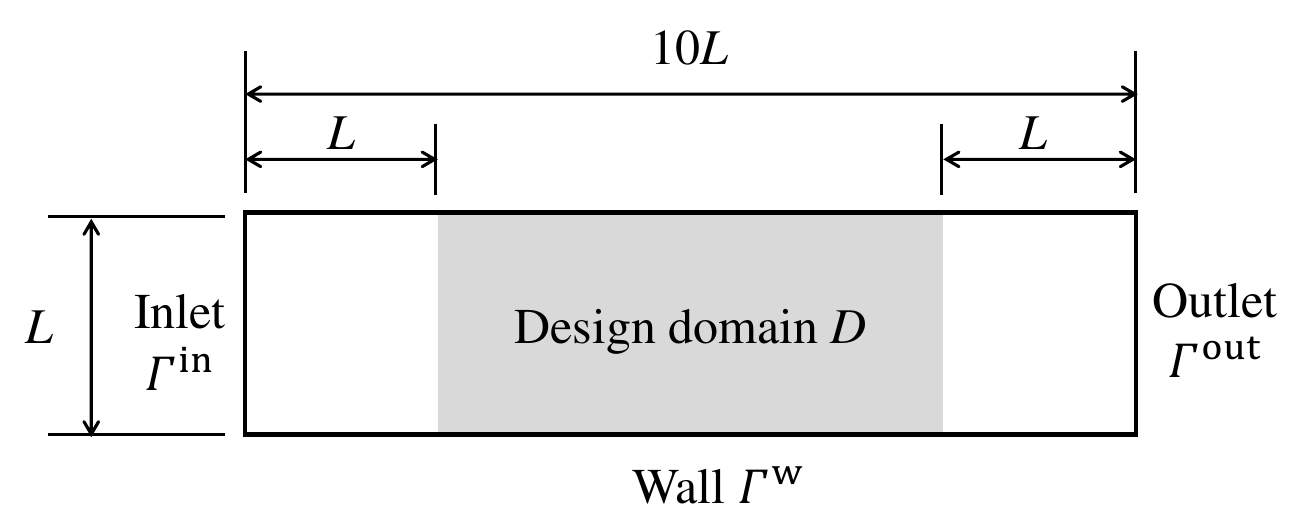} 
\caption{Dimensions of the design domain $D$ utilized for flow field design in a unit cell of the manifold microchannel heat sink.}
\label{Fig.16}
\end{figure}
\begin{figure*}[htbp]
\centering
\includegraphics[width=1.0\textwidth]{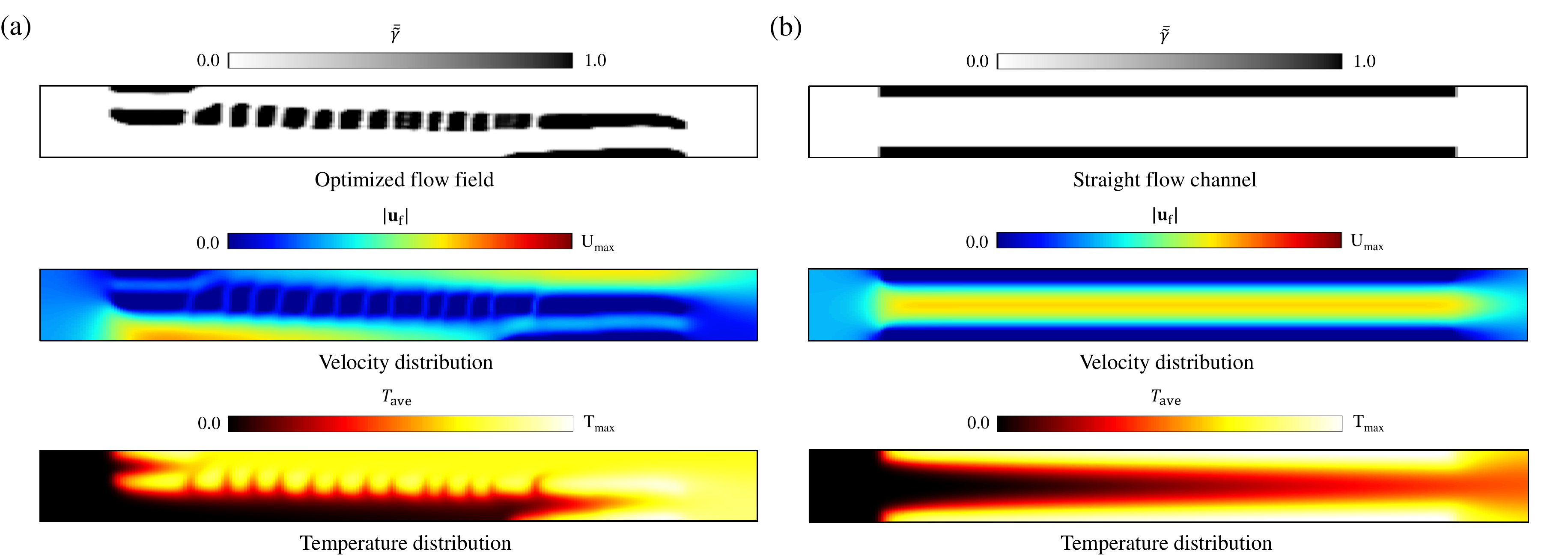} 
\caption{Comparison of the optimized flow field and straight flow channel using the nanofluid with $d_{\rm p}$ = 100 nm and $\phi_{\rm p}$ = 0.5$\%$: (a) optimized flow field (b) straight flow channel.  The maximum values on the color bar for the fluid velocity and average temperature are 0.80 m/s and 347.6 $\rm K$, respectively.}
\label{Fig.9}
\end{figure*}
\subsection{Flow field design for a manifold microchannel heat sink}
\subsubsection{Details of numerical setting}
In the second numerical example, we optimize the flow field for a unit cell of the MMCHS and compare its performance with that of a conventional parallel flow channel. Table~\ref{validation1} presents the material properties for the nanofluid, and Fig.~\ref{Fig.16} illustrates the geometric configuration of this problem, with the inlet length $L$ set to 50 $\mu$m. The analysis domain is discretized using a quadrilateral mesh of  $250 \times 30$. The pressure drop between the inlet and outlet is set to 2000 $ \rm Pa$, where the corresponding fluid Reynolds number is 70.7. The particle volume fraction at the inlet $\phi^{\rm in}_{\rm p}$ is set to 0.5$\%$. In this case, the inlet and solid temperatures are identical to those in the first numerical example.

In the optimization setup,  we set the maximum values of inverse permeability $\alpha$ and the volumetric heat generation coefficient $Q$ are $1.6 \times 10^{10}$ and $1.5 \times 10^{10}$, respectively. The initial values of design variables are set to 0.9 in the design domain. The convex parameter $q$ is set to 0.05. The filter radius $R$ is set to 4 $\rm \mu m$, and the maximum value of the projection parameter $\beta_{\rm max}$ is set to 64.
\begin{figure}[h]
\centering
\includegraphics[width=0.45\textwidth]{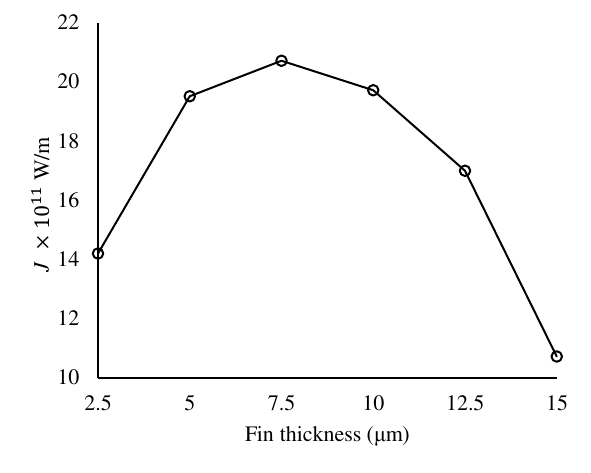} 
\caption{Objective function values of the straight flow channel with different fin thicknesses. The maximum value of $J$ is $2.07 \times 10^{11}$ W/m at a fin thickness of 7.5 $\rm \mu m$.}
\label{Fig.11}
\end{figure}
\subsubsection{Optimized flow field}
Fig.~\ref{Fig.9}(a) presents the optimization results for a unit cell of the MMCHS, with an objective function value of $2.31 \times 10^{11}$ W/m. As shown in Fig.~\ref{Fig.9}(a), the fins are primarily distributed along the central line, forming multiple branches to increase the heat dissipation area. It can be observed that the fluid flows from the inlet primarily along the bottom edge and absorbs heat as it flows around the branches. To evaluate the performance of the optimized flow field compared to a conventional parallel flow field, we constructed a straight channel within the design domain by aligning solid fins, represented by design variables, along the upper and lower boundaries. The fin thickness was adjusted to determine the optimal performance of the straight flow channel. Fig.~\ref{Fig.11} illustrates the objective function values of the straight flow channel with different fin thicknesses using the same nanofluid. The optimal performance is achieved at a fin thickness of 7.5 $\rm \mu m$, yielding an objective function value of $2.07 \times 10^{11}$ W/m. Fig.~\ref{Fig.9}(b) depicts the straight flow channel with a fin thickness of 7.5 $\rm \mu m$, along with its corresponding velocity and temperature distributions. These findings indicate that the proposed topology optimization method can generate a superior unit cell for MMCHS, achieving an 11.6$\%$ enhancement in heat transfer performance over a conventional parallel flow field.
\section{Conclusion}
\label{Conclusion}
In this study, we proposed a density-based topology optimization method for microchannel heat sinks using nanofluids. The behavior of the nanofluid was simulated using an Eulerian-Eulerian model. In the optimization setup, we formulated the problem as a heat generation maximization problem under a fixed pressure drop. The optimization was performed using the GCMMA algorithm, with sensitivities computed through automatic differentiation. In the numerical examples, we used water as the base fluid and copper as the nanoparticle material. We first validated our numerical scheme against the results from prior works. Subsequently, we examined the effects of different parameters on the optimized flow field and employed this method to design the MMCHS unit cell. The findings are summarized as follows:
\begin{itemize}
\item An increased pressure drop enhances convective effects within the flow field, resulting in more branches in the optimized flow field and enhanced heat transfer.
\item While an increase in particle volume fraction enhances the effective thermal conductivity of the nanofluid, it also increases its viscosity. Under a fixed pressure drop, the higher viscosity reduces the number of branches in the optimized flow field to maintain flow velocity. Nevertheless, the improved thermal conductivity ensures enhanced heat transfer with increasing particle volume fractions.
\item Within the scope of this study, particle diameter has a minor effect compared to particle volume fraction. Due to its minor effect on nanofluid behavior, the variations in the optimized flow field are also minimal.
\item In the MMCHS unit cell using a nanofluid with a 0.5$\%$ particle volume fraction and 100 nm particle diameter, the topology-optimized flow field achieves 11.6$\%$ higher heat transfer performance compared to a size-optimized parallel flow field.
\end{itemize}

In future work, as the simplified two-dimensional model is employed under high aspect ratio conditions in this study, we aim to focus on designing three-dimensional flow fields in microchannel heat sinks using nanofluids. This focus addresses the limitations of the current two-dimensional model and highlights the broader applicability of three-dimensional flow fields in microchannel heat sinks~\cite{fedorov2000three, ryu2003three, ho2010experimental}.
\section*{CRediT authorship contribution statement}
\textbf{Chih-Hsiang Chen:} Writing -- original draft, Conceptualization, Methodology, Software, Validation, Formal analysis, Investigation, Data curation, Visualization. \textbf{Kentaro Yaji:} Writing-- review \& editing, Conceptualization, Methodology, Supervision, Resources, Project administration, Funding acquisition.
\section*{Acknowledgement}
This work was supported by JSPS KAKENHI (Grant No. 23H01323).


\end{document}